\documentclass[12pt,reqno]{amsart}
\usepackage{enumerate, latexsym, amsmath, amsfonts, amssymb, amsthm, color}
\def\pmod #1{\ ({\rm{mod}}\ #1)}
\def\Z{\Bbb Z}

\def\l{\left}
\def\r{\right}
\def\bg{\bigg}
\def\({\bg(}
\def\){\bg)}
\def\t{\text}
\def\f{\frac}

\def\ls{\leqslant}

\def\se {\subseteq}
\def\sm{\setminus}
\def\bi{\binom}

\def\eq{\equiv}

\def\da{\delta}

\def\Proof{\noindent{\it Proof}}

\theoremstyle{plain}
\newtheorem{theorem}{Theorem}

\newtheorem{lemma}{Lemma}

\theoremstyle{definition}

\theoremstyle{remark}
\newtheorem{remark}{Remark}

\makeatletter
\@namedef{subjclassname@2020}{%
  \textup{2020} Mathematics Subject Classification}
\makeatother
 \vspace{4mm}

\begin{document}

\medskip

\title
[On some determinants arising from quadratic residues]
{On some determinants\\ arising from quadratic residues}

\author
[Chen-Kai Ren and Zhi-Wei Sun] {Chen-Kai Ren and Zhi-Wei Sun}

\address {(Chen-kai Ren) Department of Mathematics, Nanjing
University, Nanjing 210093, People's Republic of China}
\email{chRen@smail.nju.edu.cn}

\address{(Zhi-Wei Sun, corresponding author) Department of Mathematics, Nanjing
University, Nanjing 210093, People's Republic of China}
\email{zwsun@nju.edu.cn}

\keywords{Determinants, Legendre symbols, quadratic residues modulo primes.
\newline \indent 2020 {\it Mathematics Subject Classification}. Primary 11A15, 11C20; Secondary 15A15.
\newline \indent Supported by the Natural Science Foundation of China (grant 12371004).}

\begin{abstract}
Let $p>3$ be a prime, and let $d\in\Z$ with $p\nmid d$. For the determinants
$$S_m(d,p)=\det\left[(i^2+dj^2)^{m}\right]_{1\leqslant i,j \leqslant (p-1)/2}\ \ \left(\frac{p-1}2\leqslant m\leqslant p-1\right),$$
Sun recently determined $S_m(d,p)$ modulo $p$ when $m\in\{p-2,p-3\}$ and $(\frac {-d}p)=-1$.
In this paper,  we obtain $S_{p-2}(d,p)$ modulo $p$ in the remaining case $(\frac{-d}p)=1$, and determine the Legendre symbols $(\frac{S_{p-3}(d,p)}p)$ and $(\frac{S_{p-4}(d,p)}p)$ in some special cases.
\end{abstract}
\maketitle

\section{Introduction}
\setcounter{lemma}{0}
\setcounter{theorem}{0}
\setcounter{corollary}{0}
\setcounter{remark}{0}
\setcounter{equation}{0}

Let $p$ be an odd prime, and let $(\frac{.}{p})$ be the Legendre symbol.
Carlitz \cite{C} determined the characteristic polynomial of the matrix
$$\left[\left(\frac{i-j}{p}\right)\right]_{1\leqslant i,j\leqslant p-1},
$$
and Chapman \cite{C1} evaluated the determinants
$$
\det\left[\left(\frac{i+j}{p}\right)\right]_{0\leqslant i,j\leqslant (p-1)/{2}}
\ \t{and}\ \det\l[\l(\f{i+j}p\r)\r]_{1\ls i,j\ls(p-1)/2}.
$$
via quadratic Gauss sums. Vsemirnov \cite{V12,V13} confirmed a challenging conjecture of Chapman
by evaluating the determinant
$$\det\left[\left(\frac{i-j}{p}\right)\right]_{0\leqslant i,j\leqslant (p-1)/{2}}.
$$

 Let $d$ be any integer with $p\nmid d$. Sun \cite{S19} studied the determinant
$$ S(d,p)=\det\l[\l(\frac{i^2+dj^2}{p}\r)\r]_{1\leqslant i,j \leqslant (p-1)/2},
$$
and proved that $(\f{-S(d.p)}p)=1$ if $(\f dp)=1$, and $S(d,p)=0$ if $(\f{d}p)=-1$.
Grinberg, Sun and Zhao \cite{GSZ} showed that if $p>3$
then
$$ \det\l[(i^2+dj^2)\l(\frac{i^2+dj^2}{p}\r)\r]_{0\leqslant i,j \leqslant (p-1)/2}\equiv 0 \pmod p.
$$
For any integer $m$ in the interval $((p-1)/2,p-1)$, by \cite{S22} we have
$$\det[(i^2+dj^2)^m]_{0\ls i,j\ls (p-1)/2}\eq0\pmod p,$$
which extends the above Grinberg-Sun-Zhao result.

For each $m=(p-1)/2,\ldots,p-1$, Sun \cite{S24} introduced the determinant
$$S_m(d,p)=\det\left[(i^2+dj^2)^{m}\right]_{1\leqslant i,j \leqslant (p-1)/2}.$$
In 2022, Wu, She and Wang \cite{WSW} proved \cite[Conjecture 4.5]{S19} concerning the Legendre symbol
$(\frac{S_{(p+1)/2}(d,p)}{p})$.
Sun \cite{S24} proved that if $(\frac{-d}{p})=-1$ then
\begin{align*}S_{p-2}(d,p)&\eq\det\l[\f1{i^2+dj^2}\r]_{1\ls i,j\ls(p-1)/2}
\\&\equiv
\begin{cases}
  d^{(p-1)/4} \pmod p  &\text{if} \ p\equiv 1 \pmod 4,\\
   (-1)^{(p+1)/4} \pmod p &\text{if} \ p\equiv 3 \pmod 4,\\
\end{cases}
\end{align*}
and
$$S_{p-3}(d,p)\eq\det\l[\f1{(i^2+dj^2)^2}\r]_{1\ls i,j\ls(p-1)/2}
\equiv \frac{1}{4} \prod_{r=1}^{\lfloor p/4\rfloor}\l(r+\frac{1}{4}\r)^2 \pmod p.
$$

In this paper, we obtain some further results along this line. Our method is different from that of Sun \cite{S24}.

Now we present two general results.

\begin{theorem}
     Let $p>5$ be a prime, and let $d$ be an integer with $(\frac{d}{p})=-1$.
     Let $m$ be an integer in the interval $((p-1)/2,p-1)$ with $m\eq(p-1)/2\pmod2$. Then we have
 $$ S_m(d,p)\eq 0 \pmod p. $$
\end{theorem}

\begin{theorem}
    Let $p$ be a prime with $p\eq 1 \pmod 4.$ For any integer $d$ with $(\frac{d}{p})=1$ and any odd integer $m\in((p-1)/2,p-1)$, we have
    $$\l(\frac{S_{m}(d,p)}{p}\r)\neq-1.$$
\end{theorem}

Our following three theorems deal with $S_m(d,p)$ for $m=p-2,p-3,p-4$.

\begin{theorem}\label{Th1.2} Let $p$ be an odd prime, and let $d\in\Z$ with $(\frac{-d}{p})=1$. Then
\begin{equation}\label{Sp-2}
S_{p-2}(d,p)\equiv
\begin{cases}
  (-1)^{(p+3)/4}d^{(p-1)/4}(\frac{p-3}{2}!!)^2 \pmod p  &\text{if} \ p\equiv 1 \pmod 4,\\
   0 \pmod p &\text{if} \ p\equiv 3 \pmod 4.\\
\end{cases}
\end{equation}
\end{theorem}
\begin{remark} Let $p$ be a prime with $p\eq1\pmod4$, and write $p=x^2+y^2$ $(x,y\in\Z)$
with $x\eq1\pmod 4$ and $y\eq \f{p-1}2!\,x\pmod p$. (Note that $(\f{p-1}2!\,x)^2\eq-x^2\pmod p$
by Wilson's theorem.)
As $2x\eq\bi{(p-1)/2}{(p-1)/4}\pmod p$ by Gauss' congruence (cf. \cite[(9.0.1)]{BEW} or \cite{CDE}), we have
\begin{align*}2y&\eq\f{p-1}2!(2x)\eq \f{(\f{p-1}2!)^2}{(\f{p-1}4!)^2}
\\&\eq\l(2^{(p-1)/4}\f{p-3}2!!\r)^2=\l(\f 2p\r)\l(\f{p-3}2!!\r)^2\pmod p.
\end{align*}
Thus, for any $d\in\Z$ with $(\f dp)=1$, by \eqref{Sp-2} we have
$$S_{p-2}(d,p)\eq-2yd^{(p-1)/4}\pmod p.$$
With the aid of Theorem 6.2.9 of \cite[p.\,190]{BEW}, this implies Sun's conjecture (cf. \cite{S24}) that
$$S_{p-2}(1,p)\eq -2y=2\da(s,p)\sum_{k=1}^{(p-1)/2}\l(\f {k(k^2+s)}p\r)\pmod p,$$
where $s$ is any quadratic nonresidue modulo $p$, and
$$\da(s,p)=\begin{cases}1&\t{if}\ s^{(p-1)/4}\eq\f{p-1}2!\pmod p,
\\-1&\t{otherwise}.\end{cases}.$$
\end{remark}

\begin{theorem}\label{Th1.3} Let $p>3$ be a prime with $p\eq1\pmod 4$, and let $d$ be any integer with $p\nmid d $. Then
$$\l(\frac{6S_{p-3}(d,p)}{p}\r)\neq-1.$$
Moreover, if $(\frac{d}{p})=1$ and $p \equiv5 \pmod {12}$, then
$$ \l(\frac{S_{p-3}(d,p)}{p}\r)=(-1)^{(p+3)/4}.  $$
\end{theorem}
\begin{theorem}\label{Th1.4} Let $p>4$ be a prime,  and let $d$ be any integer with $(\frac{d}{p})=1$. Then
\begin{equation}\label{p-4}\l(\frac{S_{p-4}(d,p)}{p}\r)= -1 \iff
p \equiv  3,7 \pmod {20}.
 \end{equation}
\end{theorem}

We are going to provide some auxiliary results in the next section, and prove
Theorems 1.1-1.5 in Section 3.

\section{Some auxiliary results}
\setcounter{lemma}{0}
\setcounter{theorem}{0}
\setcounter{corollary}{0}
\setcounter{remark}{0}
\setcounter{equation}{0}

Let $p$ be an odd prime. By Wilson's theorem
\begin{equation}\label{(p-1)/2}(-1)^{(p+1)/2}\l(\f{p-1}2!\r)^2\eq-\prod_{k=1}^{(p-1)/2}k(p-k)=-(p-1)!\eq1\pmod p.
\end{equation}

The following lemma can be found in Sun \cite[(1.5)]{S19}.

\begin{lemma}\label{Lem2.1} For any odd prime $p$, we have
\begin{equation}\prod_{1\ls i<j\ls (p-1)/2}(j^2-i^2)\eq\begin{cases}-\f{p-1}2!\pmod p&\t{if}\ p\eq1\pmod 4,
\\1\pmod p&\t{if}\ p\eq3\pmod4.\end{cases}
\end{equation}
\end{lemma}

\begin{lemma}\label{Lem2.2} For any positive integers $m$ and $n$, we have the identity
\begin{equation}\label{ident}\prod_{1\ls i_1<\ldots<i_m\ls n}i_1\cdots i_m=(n!)^{\bi{n-1}{m-1}}.
\end{equation}
\end{lemma}
\Proof. Observe that
\begin{align*}\prod_{1\ls i_1<\ldots<i_m\ls n}i_1\cdots i_m&=\prod_{A\se\{1,\ldots,n\}\atop |A|=m}\prod_{k\in A}k
\\&=\prod_{k=1}^n k^{|\{A\se\{1,\ldots,n\}:\ k\in A\ \&\ |A\sm\{k\}|=m-1\}|}=(n!)^{\bi{n-1}{m-1}}.
\end{align*}
This proves the identity \eqref{ident}. \qed

Using the above two lemmas, we get the following auxiliary result.

\begin{theorem} Let $p$ be an odd prime. Then
\begin{equation}
\label{i^2-j^2}\prod_{1\ls i<j\ls(p-1)/2}(i^2-j^2)\l(\f1{i^2}-\f1{j^2}\r)\eq(-1)^{\lfloor p/4\rfloor}\pmod p.
\end{equation}
\end{theorem}
\Proof. Let $n=(p-1)/2$. By Lemma \ref{Lem2.2} and the congruence \eqref{(p-1)/2},
$$\prod_{1\ls i<j\ls n}(ij)^2=(n!)^{2(n-1)}\eq(-1)^{(n+1)(n-1)}=(-1)^{n-1}\pmod p.$$
By Lemma \ref{Lem2.1} and the congruence \eqref{(p-1)/2},
$$\prod_{1\ls i<j\ls n}(i^2-j^2)^2\eq(-1)^{n-1}\pmod p.$$
Therefore
\begin{align*}&\ \prod_{1\ls i<j\ls n}(i^2-j^2)\l(\f1{i^2}-\f1{j^2}\r)
\\=&\ (-1)^{\bi n2}\prod_{1\ls i<j\ls n}\f{(i^2-j^2)^2}{i^2j^2}
\eq (-1)^{n(n-1)/2}=(-1)^{\lfloor p/4\rfloor}\pmod p.
\end{align*}
This concludes the proof of \eqref{i^2-j^2}. \qed

We need the following known lemma \cite[Lemma 10]{K04} on determinants.

\begin{lemma}\label{LemP}
Let $R$ be a commutative ring with identity, and let $P(x)=\sum_{i=0}^{n-1}a_ix^i\in R[x]$. Then we have
$$\det{\left[P(X_iY_j)\right]}_{1\leqslant i,j\leqslant n}=a_0a_1\cdots a_{n-1} \prod_{1\leqslant i< j\leqslant n}(X_i-X_j)(Y_i-Y_j).$$
\end{lemma}

Now we state our second auxiliary theorem.

\begin{theorem}\label{Th1.1} Let $p=2n+1>5$ be a prime. For $d,m\in\Z$ with $p\nmid d$ and $(p-1)/2<m<p-1$, we have
$$ S_{m}(d,p) \eq a_m^2(d,p)b_m(d,p) \pmod p, $$
where
\begin{align*}
 a_m(d,p)&=\prod_{k=0}^{\lfloor (m-n-1)/2\rfloor}\l(\binom{m}{k}+\l(\frac{d}{p}\r)\binom{m}{m-n-k}\r)\\
 &\quad \times\prod_{0\leq k< n-1-\lfloor m/2 \rfloor}\binom{m}{m-n+1+k}
\end{align*}
and
$$b_m(d,p)=\begin{cases}
  (-d)^{n/2}\l(1+\l(\frac{d}{p}\r)\r)\binom{m}{(m-n)/2}\binom{m}{m/2}   &\text{if} \ 2\mid m \ and \ 2\mid n, \\
 (-1)^{(n-1)/2}\l(\frac{d}{p}\r)^{m/2}\binom{m}{m/2}  &\t{if} \ 2\mid m \ and \ 2\nmid n, \\
 (-1)^{n/2+1}d^{n/2}\l(\frac{d}{p}\r)^{(m-1)/2} &\t{if} \ 2\nmid m \ and \ 2\mid n, \\
(-1)^{(n-1)/2} \l(1+\l(\frac{d}{p}\r)\r)\binom{m}{(m-n)/2}   &\t{if} \ 2\nmid m \ and \ 2\nmid n. \\
\end{cases}
$$
\end{theorem}
\Proof. Observe that
 \begin{align*}
S_{m}(d,p)
 &=\prod_{j=1}^{n}(j^2)^{m} \times \det\l[\l(i^2j^{-2}+d\r)^{m}\r]_{1\ls i,j\ls n}\\
&=(n!)^{2m}\det\l[\l(i^2j^{-2}+d\r)^{m}\r]_{1\ls i,j\ls n}.
 \end{align*}
 By \eqref{(p-1)/2},
 $$(n!)^{2m}\eq(-1)^{m(n+1)}\pmod p.$$
For $i,j\in\{1,\ldots,n\}$, clearly
\begin{align*}
\l( i^2j^{-2}+d\r)^{m}=&\sum_{k=0}^{m}\binom{m}{k}d^{m-k}\l( i^2j^{-2}\r)^k\\
\equiv&\sum_{k=0}^{m-n}\l(\binom{m}{k}d^{m-k}+\binom{m}{n+k}d^{m-k-n}\r)\l( i^2j^{-2}\r)^k\\
&+\sum_{m-n+1\leq k< n}\binom{m}{k}d^{m-k}\l( i^2j^{-2}\r)^{k}\\
\eq & f(i^2j^{-2})\pmod p,
\end{align*}
where
\begin{align*}
 f(x)=&\sum_{k=0}^{m-n}\l(\binom{m}{k}+d^{-n}\binom{m}{m-n-k}\r)d^{m-k}x^k\\
&+\sum_{0\leq k< p-2-m}\binom{m}{m-n+1+k}d^{n-1-k}x^{m-n+1+k}.
\end{align*}
Combining the above, we obtain
\begin{equation}
\label{f3}S_{m}(d,p)\eq (-1)^{m(n+1)}\det[f(i^2j^{-2})]_{1\ls i,j\ls n}\pmod p.
\end{equation}
By Lemma \ref{LemP} and the congruence \eqref{i^2-j^2}, we have
\begin{align*}
 &\quad \det[f(i^2j^{-2})]_{1\ls i,j\ls n} \\
 &\eq (-1)^{\lfloor p/4\rfloor}\prod_{0\leq k< p-2-m}\binom{m}{m-n+1+k}d^{n-1-k}\\
 &\quad \times\prod_{k=0}^{m-n}\l(\binom{m}{k}+d^{-n}\binom{m}{m-n-k}\r)d^{m-k}\\
 &\eq (-1)^{\lfloor p/4\rfloor}d^{n(2m-n+1)/2}\prod_{0\leq k< p-2-m}\binom{m}{m-n+1+k}\\
 &\quad \times\prod_{k=0}^{m-n}\l(\binom{m}{k}+d^{-n}\binom{m}{m-n-k}\r) \pmod p.\\
\end{align*}

As $d^{-n}\equiv (\frac{d}{p})= \pm 1 \pmod p,$ we obtain
\begin{align*}
   & \l(\binom{m}{k}+d^{-n}\binom{m}{m-n-k}\r)\l(\binom{m}{m-n-k}+d^{-n}\binom{m}{k}\r)\\
 &\equiv \l(\frac{d}{p}\r) \l(\binom{m}{k}+\l(\frac{d}{p}\r)\binom{m}{m-n-k}\r)^2 \pmod p. \\
\end{align*}
When $2\mid m-n$, we have
\begin{align*}
  &\quad \prod_{k=0}^{m-n}\l(\binom{m}{k}+d^{-n}\binom{m}{m-n-k}\r)\\&\eq\prod_{k=0}^{(m-n)/2-1}\l(\binom{m}{k}+\l(\frac{d}{p}\r)\binom{m}{m-n-k}\r)^2\\
 &\quad\times \l(\frac{d}{p}\r)^{(m-n)/2} \l(1+\l(\frac{d}{p}\r)\r)\binom{m}{(m-n)/2} \pmod p  .
\end{align*}
 If $2 \nmid m-n$, then
 \begin{align*}
  &\quad \prod_{k=0}^{m-n}\l(\binom{m}{k}+d^{-n}\binom{m}{m-n-k}\r)\\
  &\eq\l(\frac{d}{p}\r)^{(m-n+1)/2}\prod_{k=0}^{(m-n-1)/2}\l(\binom{m}{k}+\l(\frac{d}{p}\r)\binom{m}{m-n-k}\r)^2
 \pmod p .  \\
\end{align*}
It is easy to verify that
$$ d^{n(2m-n+1)/2}\l(\frac{d}{p}\r)^{\lfloor (m-n+1)/2\rfloor}\eq \begin{cases} d^{n(m+1)/2} \pmod p  & \t{if} \ 2\mid m-n, \\
d^{nm/2} \pmod p  & \t{if} \ 2\nmid m-n. \\
\end{cases}
$$
Notice that
\begin{align*}&\ \prod_{0\ls k<p-2-m}\binom{m}{m-n+1+k}
\\=&\ \begin{cases}
 \bi{m}{m/2}\prod_{0\leq k< n-1-m/2}\binom{m}{m-n+1+k}^2   &\t{if} \ 2\mid m, \\
 \prod_{0\leq k< n-1-(m-1)/2}\binom{m}{m-n+1+k}^2   &\t{if} \ 2\nmid m.
\end{cases}
\end{align*}
and
$$ (-1)^{m(n+1)} (-1)^{\lfloor p/4\rfloor}=\begin{cases} (-1)^{n/2}  & \t{if} \ 2\mid m \ and \ 2\mid n, \\
(-1)^{(n-1)/2}  & \t{if} \ 2\mid m \ and \ 2\nmid n, \\
(-1)^{n/2+1}  & \t{if} \ 2\nmid m \ and \ 2\mid n, \\
(-1)^{(n-1)/2}  & \t{if} \ 2\nmid m \ and \ 2\nmid n. \\
\end{cases}
$$
Combining the above, we obtain
\begin{align*}
  S_{m}(d,p) \eq a_m^2(d,p)b_m^{'}(d,p) \pmod p,\\
\end{align*}
where
\begin{align*}
 a_m(d,p)&=\prod_{k=0}^{\lfloor (m-n-1)/2\rfloor}\l(\binom{m}{k}+\l(\frac{d}{p}\r)\binom{m}{m-n-k}\r)\\
 &\quad \times\prod_{0\leq k< n-1-\lfloor m/2 \rfloor}\binom{m}{m-n+1+k}
\end{align*}
and
$$b_m^{'}(d,p)=\begin{cases}
  (-1)^{n/2}d^{n(m+1)/2} \l(1+\l(\frac{d}{p}\r)\r)\binom{m}{(m-n)/2}\binom{m}{m/2}  & \t{if} \ 2\mid m \ and \ 2\mid n, \\
 (-1)^{(n-1)/2}d^{nm/2}\binom{m}{m/2} & \t{if} \ 2\mid m \ and \ 2\nmid n, \\
 (-1)^{n/2+1}d^{nm/2}& \t{if} \ 2\nmid m \ and \ 2\mid n, \\
 (-1)^{(n-1)/2}d^{n(m+1)/2} \l(1+\l(\frac{d}{p}\r)\r)\binom{m}{(m-n)/2}  & \t{if} \ 2\nmid m \ and \ 2\nmid n. \\
\end{cases}
$$
When $2\mid m$, we have
\begin{align*}
   &\quad(-1)^{n/2}d^{n(m+1)/2} \l(1+\l(\frac{d}{p}\r)\r)\\
   &\eq (-1)^{n/2}d^{n/2}\l(\frac{d}{p}\r)^{m/2}\l(1+\l(\frac{d}{p}\r)\r)\\
   &\eq(-d)^{n/2}\l(1+\l(\frac{d}{p}\r)\r) \pmod p \\
\end{align*}
and
$$(-1)^{(n-1)/2}d^{nm/2}\eq(-1)^{(n-1)/2}\l(\frac{d}{p}\r)^{m/2}\pmod p.
$$
If $2\nmid m$, then
\begin{align*}
     &\quad(-1)^{(n-1)/2}d^{n(m+1)/2} \l(1+\l(\frac{d}{p}\r)\r) \\
     &\eq (-1)^{(n-1)/2}\l(\frac{d}{p}\r)^{(m+1)/2} \l(1+\l(\frac{d}{p}\r)\r)\\
     &\eq  (-1)^{(n-1)/2} \l(1+\l(\frac{d}{p}\r)\r) \pmod p\\
\end{align*}
and
$$(-1)^{n/2+1}d^{nm/2}\eq (-1)^{n/2+1}d^{n/2}\l(\frac{d}{p}\r)^{(m-1)/2} \pmod p .
$$
Thus
$$b_m^{'}(d,p)\eq b_m(d,p) \pmod p$$
where
$$b_m(d,p)=\begin{cases}
  (-d)^{n/2}\l(1+\l(\frac{d}{p}\r)\r)\binom{m}{(m-n)/2}\binom{m}{m/2}  & \t{if} \ 2\mid m \ and \ 2\mid n, \\
 (-1)^{(n-1)/2}\l(\frac{d}{p}\r)^{m/2}\binom{m}{m/2} & \t{if} \ 2\mid m \ and \ 2\nmid n, \\
 (-1)^{n/2+1}d^{n/2}\l(\frac{d}{p}\r)^{(m-1)/2}& \t{if} \ 2\nmid m \ and \ 2\mid n, \\
(-1)^{(n-1)/2} \l(1+\l(\frac{d}{p}\r)\r)\binom{m}{(m-n)/2}  & \t{if} \ 2\nmid m \ and \ 2\nmid n. \\
\end{cases}
$$
Therefore,
$$ S_{m}(d,p) \eq a_m^2(d,p)b_m(d,p) \pmod p.$$ \qed

\section{Proofs of Theorems 1.1-1.5}
\setcounter{lemma}{0}
\setcounter{theorem}{0}
\setcounter{corollary}{0}
\setcounter{remark}{0}
\setcounter{equation}{0}

Theorem 1.1 follows immediately from Theorem 2.2.

\medskip
\noindent{\bf Proof of Theorem 1.2}.
Since $(\frac{-1}{p})=1$, by Theorem 2.2 we obtain
$$\l(\frac{b_m(d,p)}{p}\r)=\l(\frac{d}{p}\r)^{n/2}=1 .$$
Thus
$$\l(\frac{S_m(d,p)}{p}\r)=\l(\frac{a_m(d,p)}{p}\r)^2\l(\frac{b_m(d,p)}{p}\r)=\l(\frac{a_m(d,p)}{p}\r)^2\neq -1 .$$ \qed

\medskip
\noindent{\bf Proof of Theorem 1.3}. Set $n=(p-1)/2$.  When $p\eq 3 \pmod 4 $ and $(\frac{d}{p})=-1$, we have
$$ S_{p-2}(d,p)\eq 0 \pmod p.$$
by Theorem 1.1.

Now we assume $p\eq 1 \pmod 4$. Then $(\frac{d}{p})=1$ and $2\mid n$. By Theorem 2.2, we have
 $$ S_{p-2}(d,p) \eq a_{p-2}^2(d,p)b_{p-2}(d,p) \pmod p,$$
where
\begin{align*}
 a_{p-2}(d,p)&=\prod_{k=0}^{n/2-1}\l(\binom{p-2}{k}+\binom{p-2}{n-1-k}\r)\\
\end{align*}
and
$$ b_{p-2}(d,p)=(-1)^{n/2+1}d^{n/2}. $$
Since
    $$ \binom{p-r-1}{k}\equiv \bi{-r-1}{k}\eq(-1)^k\binom{k+r}{r} \pmod p,
    $$
we can verify that
\begin{align*}
  \binom{p-2}{k}+\binom{p-2}{n-1-k}&\equiv (-1)^k\binom{k+1}{1}+(-1)^{n-1-k} \binom{n-k}{1} \\
  &\equiv  (-1)^k\l( 2k-n+1 \r) \pmod p.
  \end{align*}
Combining the above, we obtain
\begin{align*}
 S_{p-2}(d,p)&\eq  (-1)^{n/2+1}d^{n/2}\prod_{k=0}^{n/2-1}\l(\binom{p-2}{k}+d^{-n}\binom{p-2}{n-1-k}\r)^2\\
& \equiv (-1)^{n/2+1}d^{n/2}\prod_{k=0}^{n/2-1} \l((-1)^k( 2k-n+1)\r)^2  \\
  & \equiv  (-1)^{(p+3)/4}d^{(p-1)/4}\l(\l(\frac{p-3}{2}\r)!!\r)^2 \pmod{p}.
\end{align*}
This concludes the proof. \qed

\medskip
\noindent{\bf Proof of Theorem 1.4}.
Set $n=(p-1)/2$. By \cite[Theorem 1.2]{S24}, if $(\frac{d}{p})=-1$ then
 $$\l(\frac{S_{p-3}(d,p)}{p}\r)=0.$$
When $(\frac{d}{p})=1,$ by Theorem 2.2 we obtain
$$ S_{p-3}(d,p) \eq a_{p-3}^2(d,p)b_{p-3}(d,p) \pmod p,$$
where
\begin{align*}
 a_{p-3}(d,p)&=\prod_{k=0}^{n/2-2}\l(\binom{p-3}{k}+\binom{p-3}{n-2-k}\r)\\
\end{align*}
and
$$ b_{p-3}(d,p)=2(-d)^{n/2}\binom{p-3}{(n-2)/2}\binom{p-3}{n-1}. $$

Observe that
\begin{align*}
  2\binom{p-3}{(n-2)/2}\binom{p-3}{n-1} &\equiv 2(-1)^{n/2-1}\bi{(n-2)/2}{2}(-1)^{n-1}\binom{n-1}{2} \\
  &\eq\frac{ (-1)^n675}{128} \pmod p.
\end{align*}
 If there exists an integer $k$ with $0\leq k\leq n/2-2$ satisfying $$\binom{p-3}{k}+\binom{p-3}{n-2-k} \equiv 0 \pmod p, $$
 then
 $$\l(\frac{6S_{p-3}(d,p)}{p}\r)=0 \neq-1.$$
 Otherwise, combining the above, we obtain
\begin{align*}
 \l(\frac{6S_{p-3}(d,p)}{p}\r) &=\l(\frac{d^{n/2}}{p}\r)=\l(\frac{d}{p}\r)^{n/2} =1.
\end{align*}
We conclude that if $p\equiv 1 \pmod 4$, then
$$\l(\frac{6S_{p-3}(d,p)}{p}\r)\neq-1.$$

Suppose that $p\equiv 5 \pmod {12}$ and $(\frac{d}{p})=1$. We claim that
$$p \nmid \l( \binom{p-3}{k}+\binom{p-3}{n-2-k}\r) $$
 for any integer $k$  with $0\ls k\ls n/2-2.$
Note that
\begin{align*}
  \binom{p-3}{k}+\binom{p-3}{n-2-k}
 &\equiv (-1)^k\binom{k+2}{2}+(-1)^{n-2-k} \binom{n-k}{2} \\
  &\equiv \frac{(-1)^k}{2}\l( (k+1)(k+2)+(k-n)(k-n-1)   \r)\\
  &\equiv  (-1)^k\l( \l(k+\frac{5}{4}\r)^2-\frac{3}{16} \r) \pmod p.
  \end{align*}
Since $(\frac{3}{p})=-1$, the claim holds.

 From the above discussion, we obtain
\begin{align*}
 \l(\frac{S_{p-3}(d,p)}{p}\r)
 =\l(\frac{6}{p}\r)
 =(-1)^{\frac{p+3}{4}}.
\end{align*}
This ends the proof. \qed

\medskip
\noindent{\bf Proof of Theorem 1.5}. Set $n=(p-1)/2$.
If
$$\l(\frac{S_{p-4}(d,p)}{p}\r)= -1,$$ then $p\equiv 3 \pmod 4$ by Theorem 1.2.
Thus, by Theorem 2.2 we obtain
$$ S_{p-4}(d,p) \eq a_{p-4}^2(d,p)b_{p-4}(d,p) \pmod p,$$
where
\begin{align*}
 a_{p-4}(d,p)&=\binom{p-4}{n-2}\prod_{k=0}^{(n-5)/2}\l(\binom{p-4}{k}+\binom{p-4}{n-3-k}\r)\\
\end{align*}
and
$$ b_{p-4}(d,p)=2(-1)^{(n-1)/2} \binom{p-4}{(n-3)/2}. $$
It is easy to verify that
\begin{align*}
    2\binom{p-4}{(n-3)/2}&\eq 2(-1)^{(n-3)/2} \binom{(n+3)/2}{3} \\
    & \eq \frac{5(-1)^{(n-1)/2}}{64} \pmod p.
\end{align*}
If there exists an integer $k$ with $0\leq k\leq (n-5)/2$ satisfying $$\binom{p-4}{k}+\binom{p-4}{n-3-k} \equiv 0 \pmod p, $$
 then
 $$\l(\frac{S_{p-4}(d,p)}{p}\r)=0 \neq-1.$$
 Otherwise, combining the above, we obtain
\begin{align*}
\l(\frac{S_{p-4}(d,p)}{p}\r)=\l(\frac{5}{p}\r).
 \end{align*}
 Thus, we deduce that $p\eq \pm 2 \pmod 5.$ This proves the ``$\Rightarrow$" direction of \eqref{p-4}.

Now suppose that $p\equiv 3 \pmod 4$ and $p\equiv \pm 2 \pmod 5$. We claim that
 $$p \nmid \l( \binom{p-4}{k}+\binom{p-4}{n-3-k}\r) $$
 for any
 integer $k$  with $0\ls k\ls (n-5)/2.$ Note that
\begin{align*}
    &\binom{p-4}{k}+\binom{p-4}{n-3-k}\\
    \equiv& (-1)^k\binom{k+3}{3}+(-1)^{n-3-k} \binom{n-k}{3} \\
  \equiv& \frac{(-1)^k}{6}\l( \l(k+1\r)\l(k+2\r)\l(k+3\r)-\l(k+\frac{1}{2}\r)\l(k+\frac{3}{2}\r)\l(k+\frac{5}{2}\r)   \r)\\
  \equiv&  \frac{(-1)^k}{4}\l( \l(k+\frac{7}{4}\r)^2-\frac{5}{16} \r) \pmod p.
\end{align*}
Since $(\frac{5}{p})=-1$, the claim holds.
In view of this, we obtain
 $$\l(\frac{S_{p-4}(d,p)}{p}\r)=\l(\frac{5}{p}\r)=-1.$$
 So the ``$\Leftarrow$" direction of \eqref{p-4} also holds.\qed

\setcounter{conjecture}{0}
\end{document}